\documentclass[12pt]{article} 
\usepackage{amsfonts,amsmath,latexsym,amssymb,mathrsfs,amsthm,comment,setspace}
\usepackage{slashbox}
\usepackage{caption}

\let\OLDthebibliography\thebibliography
\renewcommand\thebibliography[1]{
  \OLDthebibliography{#1}
  \setlength{\parskip}{0pt}
  \setlength{\itemsep}{0pt plus 0.0ex}
}

\def\numberlikeadb{\global\def\theequation{\thesection.\arabic{equation}}}
\numberlikeadb
\newtheorem{theorem}{Theorem}[section]
\newtheorem{lemma}[theorem]{Lemma}
\newtheorem{corollary}[theorem]{Corollary}

\newtheorem{proposition}[theorem]{Proposition}
\newtheorem{remark}[theorem]{Remark}

\usepackage{color}

\usepackage{lscape}
\usepackage{caption}
\usepackage{multirow}
\begin{document}

\title{Bounds for an integral involving the modified Lommel function of the first kind}
\author{Robert E. Gaunt\footnote{Department of Mathematics, The University of Manchester, Oxford Road, Manchester M13 9PL, UK, robert.gaunt@manchester.ac.uk}}

\date{\today} 
\maketitle

\vspace{-5mm}


\begin{abstract}Simple upper and lower bounds are established for the integral $\int_0^x\mathrm{e}^{-\beta u}u^\nu t_{\mu,\nu}(u)\,\mathrm{d}u$, where $x>0$, $0<\beta<1$, $\mu+\nu>-2$, $\mu-\nu\geq-3$, and $t_{\mu,\nu}(x)$ is the modified Lommel function of the first kind. Our bounds complement and improve on existing bounds for this integral, by either being sharper or increasing the range of validity. Our bounds also generalise recent bounds for an integral involving the modified Struve function of the first kind, and in some cases more a direct approach lead to sharper bounds when our general bounds are specialised to the modified Struve case.
\end{abstract}


\noindent{{{\bf{AMS 2010 Subject Classification:}}} Primary 33C20; 26D15

\section{Introduction}\label{intro}

In the recent papers \cite{gaunt ineq1,gaunt ineq3,gaunt ineq8}, simple upper and lower bounds, involving the modified Bessel function of the first kind $I_\nu(x)$, were obtained for the integral
\begin{equation}\label{intbes}\int_0^x \mathrm{e}^{-\beta u} u^\nu I_\nu(u)\,\mathrm{d}u, 
\end{equation}
where $x>0$ and $0\leq\beta<1$.  The conditions imposed on $\nu$ differed from one inequality to another, but, in order to ensure that the integral exists, it was always assumed that $\nu>-\frac{1}{2}$.  For $0<\beta<1$, simple closed-form formulas are not available for this integral. The bounds of \cite{gaunt ineq1,gaunt ineq3,gaunt ineq8} were required in the development of Stein's method \cite{chen,np12,stein} for variance-gamma approximation \cite{
eichelsbacher, gaunt vg, gaunt vg2,gaunt vg3}. Despite their simple form, the bounds of \cite{gaunt ineq1,gaunt ineq3,gaunt ineq8} are quite accurate, and so may also be effective in other problems concerning modified Bessel functions; see, for example, \cite{bs09,baricz3} in which inequalities for modified Bessel functions of the first kind were used to establish tight bounds for the generalized Marcum $Q$-function, encountered in radar signal processing.

 A natural analogue of the problem considered by \cite{gaunt ineq1,gaunt ineq3,gaunt ineq8} is to ask for simple bounds, involving the modified Lommel function of the first kind $t_{\mu,\nu}(x)$, for the integral
\begin{equation}\label{intstruve}\int_0^x \mathrm{e}^{-\beta u} u^{\nu} t_{\mu,\nu}(u)\,\mathrm{d}u, 
\end{equation}
where $x>0$, $0\leq\beta<1$, $\mu+\nu>-2$, $\mu-\nu\geq-3$ (that is $\mu>-\frac{5}{2}$, $-\mu-2<\nu\leq\mu+3$). The condition $\mu+\nu>-2$ assures us that the integral exists, and the conditions $\mu+\nu>-2$, $\mu-\nu\geq-3$ together ensure that the integrand (when normalised according to (\ref{tildedefn})) is positive.  This problem was recently studied by \cite{gaunt lommel1}, and will also be the focus of this paper.
 Up to a multiplicative constant, the modified Lommel function $t_{\mu,\nu}(x)$ generalises the modified Struve function of the first kind $\mathbf{L}_\nu(x)$ (see (\ref{bnm})). Therefore, the problem of asking for simple bounds for the integral (\ref{intstruve}) is also a natural generalisation of the problem of obtaining simple bounds, in terms of the modified Struve function of the first kind, for the integral
 \begin{equation}\label{intstruve0}\int_0^x \mathrm{e}^{-\beta u} u^{\nu} \mathbf{L}_\nu(u)\,\mathrm{d}u, 
\end{equation}
where $x>0$, $0\leq\beta<1$, $\nu>-1$, that was recently studied by \cite{gaunt ineq4,gaunt struve int2}. 
 

The modified Lommel function of the first kind $t_{\mu,\nu}(x)$ is defined by the hypergeometric series
\begin{equation*}t_{\mu,\nu}(x)=2^{\mu-1}\Gamma\big(\tfrac{\mu-\nu+1}{2}\big)\Gamma\big(\tfrac{\mu+\nu+1}{2}\big)\sum_{k=0}^\infty\frac{(\frac{1}{2}x)^{\mu+2k+1}}{\Gamma\big(k+\frac{\mu-\nu+3}{2}\big)\Gamma\big(k+\frac{\mu+\nu+3}{2}\big)}.
\end{equation*}
It arises as a particular solution of the modified Lommel differential equation \cite{r64,s36}
\begin{equation*}x^2f''(x)+xf'(x)-(x^2+\nu^2)f(x)=x^{\mu+1}.
\end{equation*}
In the literature, different notation can be found for modified Lommel functions; we use that of \cite{zs13}.  The terminology modified Lommel function of the \emph{first kind} is also not yet standard in the literature, but was recently introduced by \cite{gaunt lommel}. This is consistent with the terminology Lommel function of the \emph{first kind} used by \cite{bk14} for the function $s_{\mu,\nu}(x)$, which is related to  modified Lommel function of the first kind through $t_{\mu,\nu}(x)=-\mathrm{i}^{1-\mu}s_{\mu,\nu}(\mathrm{i}x)$ (see \cite{r64,zs13}). Modified Lommel functions are encountered in areas of the applied sciences as diverse as the theory of steady-state temperature distribution \cite{g29}, stress distributions in cylindrical objects \cite{s85} and scattering amplitudes in quantum optics \cite{t73}; for a list of further applications see \cite{gaunt lommel}. The modified Struve function of the first kind $\mathbf{L}_{\nu}(x)$, which is an important special case of the modified Lommel function $t_{\mu,\nu}(x)$, is itself a widely used special functions; see \cite{bp14} for a list of some of its application areas.

It will be convenient to follow \cite{gaunt lommel} and use the following normalisation for the purpose of removing cumbersome multiplicative constants from our calculations:
\begin{equation}\label{tildedefn}\tilde{t}_{\mu,\nu}(x)=\frac{t_{\mu,\nu}(x)}{2^{\mu-1}\Gamma\big(\frac{\mu-\nu+1}{2}\big)\Gamma\big(\frac{\mu+\nu+1}{2}\big)}=\sum_{k=0}^\infty\frac{(\frac{1}{2}x)^{\mu+2k+1}}{\Gamma\big(k+\frac{\mu-\nu+3}{2}\big)\Gamma\big(k+\frac{\mu+\nu+3}{2}\big)}.
\end{equation}
To ease the exposition, we will also refer to $\tilde{t}_{\mu,\nu}(x)$ as the modified Lommel function of the first kind.  From now on, we will work with the function $\tilde{t}_{\mu,\nu}(x)$, with   results for $t_{\mu,\nu}(x)$ being easily deduced. For example, in interpreting the bounds obtained in this paper, it is useful to note that 
\begin{align*}\frac{t_{\mu+1,\nu+1}(x)}{\tilde{t}_{\mu+1,\nu+1}(x)}=(\mu+\nu+1)\frac{t_{\mu,\nu}(x)}{\tilde{t}_{\mu,\nu}(x)}.
\end{align*}
We also note the important special case
\begin{equation}\label{bnm}\tilde{t}_{\nu,\nu}(x)= \mathbf{L}_\nu(x).
\end{equation}
A list of further basic properties of the modified Lommel function $\tilde{t}_{\mu,\nu}(x)$ that are needed in this paper is given in Appendix \ref{appa}. 



When $\beta=1$, the integral (\ref{intstruve}) can be evaluated exactly in terms of the modified lommel function $\tilde{t}_{\mu,\nu}(x)$. Suppose $\mu>-\frac{3}{2}$, $-\frac{1}{2}<\nu<\mu+1$. Then, for $x>0$,
\begin{align}\int_0^x \mathrm{e}^{-u}u^\nu \tilde{t}_{\mu,\nu}(u)\,\mathrm{d}u&=\frac{\mathrm{e}^{-x}x^{\nu+1}}{2\nu+1}\big(\tilde{t}_{\mu,\nu}(x)+\tilde{t}_{\mu+1,\nu+1}(x)\big)\nonumber\\
\label{intfor}&\quad-\frac{\gamma(\mu+\nu+2,x)}{2^\mu(2\nu+1)\Gamma\big(\frac{\mu-\nu+1}{2}\big)\Gamma\big(\frac{\mu+\nu+3}{2}\big)},
\end{align}
where  $\gamma(a,x)=\int_0^x \mathrm{e}^{-u}u^{a-1}\,\mathrm{d}u$ is the lower incomplete gamma function. We will make use of this integral formula later in proving one of the main results of this paper. A simple verification of the formula (\ref{intfor}) is given in Appendix \ref{appb}. The condition $\mu>-\frac{3}{2}$, $-\frac{1}{2}<\nu<\mu+1$ ensures that $\mu+\nu>-2$ (so that the integral exists), $\nu>-\frac{1}{2}$ and $\mu-\nu>-1$. When $0\leq\beta<1$, an exact formula for the integral (\ref{intstruve}) in terms of the function $\tilde{t}_{\mu,\nu}(x)$ is not available, although the integral can be evaluated exactly in terms of the generalized hypergeometric function (the case $\beta=0$) or an infinite series involving lower incomplete gamma functions (when $0<\beta<1$); see \cite{gaunt lommel1}. 
The fact that no simple closed-form formulas are available for the integral (\ref{intstruve}) provides the 
motivation for establishing simple bounds, involving the function $ \tilde{t}_{\mu,\nu}(x)$ itself, for this integral.

Several upper bounds and a lower bound for the integral (\ref{intstruve}) were established by \cite{gaunt lommel1}.  In this paper, we complement that work by obtaining new bounds for the integral (\ref{intstruve}) that are either sharper or have a larger range of validity than those of \cite{gaunt lommel1}. We achieve our bounds by adapting the approach used in the recent paper \cite{gaunt struve int2} in which similar improvements were obtained on the bounds of \cite{gaunt ineq4} for the integral (\ref{intstruve0}) involving $\mathbf{L}_\nu(x)$. This approach is fruitful because the properties of the modified Struve function $\mathbf{L}_\nu(x)$ used by \cite{gaunt struve int2} generalise naturally to the modified Lommel function $\tilde{t}_{\mu,\nu}(x)$. As such, our bounds generalise those of \cite{gaunt struve int2}, although in some cases, by using more direct arguments, some of our bounds when specialised to the modified Struve function $\mathbf{L}_\nu(x)$ are in fact sharper than those of \cite{gaunt struve int2}.

In Theorem \ref{tiger2}, we extend the range of validity (from $\mu>-\frac{1}{2}$, $\frac{1}{2}\leq\nu<\mu+1$ to $\mu>-\frac{3}{2}$, $-\frac{1}{2}<\nu<\mu+1$) of the upper bounds of \cite{gaunt lommel1} for the integral (\ref{intstruve}), with our bounds having the same functional form, but larger numerical constants. When specialised to the case of the modified Struve function $\mathbf{L}_\nu(x)$, our bounds are sharper than those of \cite{gaunt struve int2} (see Corollary \ref{struvecor}). 
We also establish several lower bounds for the integral (\ref{intstruve}) (Theorem \ref{tiger1}), one of which is strictly sharper than the only lower bound given in \cite{gaunt lommel1}.  In fact, all lower bounds derived in this paper are tight in the limit $x\rightarrow\infty$, a property not enjoyed by the lower bound of \cite{gaunt lommel1}.  When specialised to the case of the modified Struve function $\mathbf{L}_\nu(x)$, one of the lower bounds is sharper than that of \cite{gaunt struve int2} (see Corollary \ref{struvecor}). 
All our results are stated in Section \ref{sec2}. A discussion and comparison of the bounds of Section \ref{sec2} is given in Section \ref{secnew}. The proofs of the results of Section \ref{sec2} are given in Sections \ref{sec3}. Some basic properties of the modified Lommel function $\tilde{t}_{\mu,\nu}(x)$ that are needed in the paper are collected in Appendix \ref{appa}. Finally, in Appendix \ref{appb}, we provide a simple verification of the integral formula (\ref{intfor}).

\section{Main results}\label{sec2}

We begin by introducing some notation, and a lemma (which, as with all other results stated in this section, is proved in Section \ref{sec3}) that is used in the proof of Theorem \ref{tiger2}. Firstly, we let
\begin{align*}A_{\mu,\nu}&=\frac{1}{2}(\mu-\nu)+\frac{1}{2}\sqrt{(\mu-\nu)^2+8(\mu+\nu+3)}, \\
B_{\mu,\nu,\beta}&=\frac{\beta^{-\mu-\nu-1}\gamma(\mu+\nu+1,\beta x)}{2^\mu\Gamma\big(\frac{\mu-\nu+1}{2}\big)\Gamma\big(\frac{\mu+\nu+3}{2}\big)},
\end{align*} 
where $\gamma(a,x)=\int_0^x\mathrm{e}^{-u}u^{a-1}\,\mathrm{d}u$ is  the lower incomplete gamma function.

\begin{lemma}\label{lem1}Let $\mu>-\frac{3}{2}$, $-\frac{1}{2}<\nu<\mu+1$ and $0<\beta<1$. Fix $x_*>0$. Then, for $0<x\leq x_*$,
\begin{equation}\label{ineqb1ccc}\int_0^x \mathrm{e}^{-\beta u}u^{\nu}\tilde{t}_{\mu,\nu}(u)\,\mathrm{d}u<\frac{\mu+\nu+3+2x_*}{2\nu+1}\mathrm{e}^{-\beta x} x^\nu \tilde{t}_{\mu+1,\nu+1}(x).
\end{equation}
Now, fix $x_*>\frac{1}{1-\beta}$.  Then, for $x\geq x_*$,
\begin{equation}\label{ineqb1}\int_0^x \mathrm{e}^{-\beta u}u^{\nu}\tilde{t}_{\mu,\nu}(u)\,\mathrm{d}u<M_{\mu,\nu,\beta}(x_*)\mathrm{e}^{-\beta x} x^\nu \tilde{t}_{\mu+1,\nu+1}(x), 
\end{equation}
where
\begin{equation}\label{mng}M_{\mu,\nu,\beta}(x_*)=\max\bigg\{\frac{\mu+\nu+3+2x_*}{2\nu+1},\frac{x_*}{(1-\beta)x_*-1}\bigg\}.
\end{equation}
\end{lemma}

The bounds given in Theorems \ref{tiger2} and \ref{tiger1} below complement and improve on bounds of \cite{gaunt lommel1} for the integral (\ref{intstruve}). When specialised to the case of the modified Struve function $\mathbf{L}_\nu(x)$, the bounds also generalise, and in some cases are sharper than, bounds of \cite{gaunt struve int2} for the integral (\ref{intstruve0}) involving the function $\mathbf{L}_\nu(x)$. The bounds of  Theorems \ref{tiger2} and \ref{tiger1} are also natural analogues of bounds recently established by \cite{gaunt ineq8} for the related integral $\int_0^x \mathrm{e}^{-\beta u}u^{\nu}I_\nu(u)\,\mathrm{d}u$. 




\begin{theorem}\label{tiger2}Let $\mu>-\frac{3}{2}$, $-\frac{1}{2}<\nu<\mu+1$ and $0<\beta<1$. 
Then, for $x>0$,
\begin{align}\label{ineqb1ss}\int_0^x\mathrm{e}^{-\beta u}u^\nu \tilde{t}_{\mu,\nu}(u)\,\mathrm{d}u&<\frac{\mu+\nu+3+A_{\mu,\nu}}{(2\nu+1)(1-\beta)}\mathrm{e}^{-\beta x} x^\nu \tilde{t}_{\mu+1,\nu+1}(x), \\
\label{ineqb10}\int_0^x\mathrm{e}^{-\beta u}u^\nu \tilde{t}_{\mu,\nu}(u)\,\mathrm{d}u&<\frac{\mu+\nu+\frac{1}{2}(9+\sqrt{17})}{(2\nu+1)(1-\beta)}\mathrm{e}^{-\beta x}x^\nu \tilde{t}_{\mu+1,\nu+1}(x).
\end{align}
\end{theorem}



\begin{theorem}\label{tiger1}Let $0<\beta<1$. Then, for $x>0$,
\begin{align}\label{ineqb2}\int_0^x \mathrm{e}^{-\beta u}u^{\nu}\tilde{t}_{\mu,\nu}(u)\,\mathrm{d}u&>\frac{1}{1-\beta}\big\{\mathrm{e}^{-\beta x}x^\nu\tilde{t}_{\mu,\nu}(x)-B_{\mu,\nu,\beta}\big\},   \\
\label{ineqb3}\int_0^x \mathrm{e}^{-\beta u}u^{\nu}\tilde{t}_{\mu,\nu}(u)\,\mathrm{d}u&>\frac{1}{1-\beta}\bigg\{\bigg(1-\frac{4\nu^2}{(2\nu-1)(1-\beta)}\frac{1}{x}\bigg)\mathrm{e}^{-\beta x}x^\nu \tilde{t}_{\mu,\nu}(x)-B_{\mu,\nu,\beta}\bigg\},\\
\int_0^x\mathrm{e}^{-\beta u}u^\nu \tilde{t}_{\mu,\nu}(u)\,\mathrm{d}u&>\frac{1}{1-\beta}\bigg\{\bigg(1-\frac{2\nu(\mu+\nu+\frac{1}{2}(5+\sqrt{17}))}{(2\nu-1)(1-\beta)}\frac{1}{x}\bigg)\times\nonumber\\
\label{ineqb12}&\quad\times\mathrm{e}^{-\beta x}x^\nu \tilde{t}_{\mu,\nu}(x)-B_{\mu,\nu,\beta}\bigg\}, \\
\label{ineqb4}\int_0^x \mathrm{e}^{-\beta u}u^{\nu}\tilde{t}_{\mu,\nu}(u)\,\mathrm{d}u&>\mathrm{e}^{-\beta x}x^\nu\sum_{k=0}^\infty \beta^k \tilde{t}_{\mu+k+1,\nu+k+1}(x),
\end{align}
Inequality (\ref{ineqb2}) is valid for $\mu>-1$, $-\mu-1<\nu\leq0$; inequality (\ref{ineqb3}) is valid for $\mu>\frac{1}{2}$, $\frac{3}{2}\leq\nu<\mu+1$; inequality (\ref{ineqb12}) is valid for $\mu>-\frac{1}{2}$, $\frac{1}{2}<\nu<\mu+1$; and inequality (\ref{ineqb4}) holds for $\mu>-\frac{5}{2}$, $-\mu-2<\nu\leq\mu+3$.  Inequalities (\ref{ineqb2})--(\ref{ineqb4}) are tight in the limit $x\rightarrow\infty$. 
\end{theorem}

\begin{remark}The series in the lower bound in (\ref{ineqb4}) can be expressed in terms of the Kamp\'e de F\'eriet function (defined in Appendix \ref{appa}) so that the lower bound is equal to
\begin{align*}&\mathrm{e}^{-\beta x}x^\nu\sum_{k=0}^\infty \beta^k \tilde{t}_{\mu+k+1,\nu+k+1}(x)\\
&\quad= \frac{\mathrm{e}^{-\beta x}x^{\mu+\nu+2} }{2^{\mu+2} \Gamma(\frac{\mu-\nu+3}{2}) \Gamma(\frac{\mu+\nu+5}{2})} \mathrm{F}_{1:-;1}^{-:1;1} \Bigg[ \begin{array}{c} -:1;1\\ \frac{\mu+\nu+5}{2}:-; \frac{\mu-\nu+3}{2}
   \end{array} \Bigg| \frac{\beta x}2, \frac{x^2}4 \Bigg]\,. 
   \end{align*}
   \end{remark}



The bounds in the following proposition are sharper than the bounds (\ref{ineqb2})--(\ref{ineqb12}), since $\tilde{t}_{\mu+1,\nu+1}(x)<\tilde{t}_{\mu,\nu}(x)$, $x>0$, $\mu>-\frac{3}{2}$, $-\frac{1}{2}\leq\nu<\mu+1$ (see (\ref{Imon})).

\begin{proposition}\label{propone}Let $0<\beta<1$. Then, for $x>0$,
\begin{align}\label{ineqb21}\int_0^x \mathrm{e}^{-\beta u}u^{\nu}\tilde{t}_{\mu+1,\nu+1}(u)\,\mathrm{d}u&>\frac{1}{1-\beta}\big\{\mathrm{e}^{-\beta x}x^\nu\tilde{t}_{\mu,\nu}(x)-B_{\mu,\nu,\beta}\big\}, \\
\int_0^x \mathrm{e}^{-\beta u}u^{\nu}\tilde{t}_{\mu+1,\nu+1}(u)\,\mathrm{d}u&>\frac{1}{1-\beta}\bigg\{\bigg(1-\frac{4\nu^2}{(2\nu-1)(1-\beta)}\frac{1}{x}\bigg)\mathrm{e}^{-\beta x}x^\nu \tilde{t}_{\mu,\nu}(x)\nonumber\\
\label{ineqb22}&\quad-B_{\mu,\nu,\beta}\bigg\}, \\
\int_0^x\mathrm{e}^{-\beta u}u^\nu \tilde{t}_{\mu+1,\nu+1}(u)\,\mathrm{d}u&>\frac{1}{1-\beta}\bigg\{\bigg(1-\frac{2\nu(\mu+\nu+\frac{1}{2}(5+\sqrt{17}))}{(2\nu-1)(1-\beta)}\frac{1}{x}\bigg)\times\nonumber\\
\label{ineqb23}&\quad\times\mathrm{e}^{-\beta x}x^\nu \tilde{t}_{\mu,\nu}(x) -B_{\mu,\nu,\beta}\bigg\}.
\end{align}
Inequality (\ref{ineqb21}) is valid for $\mu>-1$, $-\mu-1<\nu\leq0$; inequality (\ref{ineqb22}) is valid for $\mu>\frac{1}{2}$, $\frac{3}{2}\leq\nu<\mu+1$; and inequality (\ref{ineqb23}) holds for $\mu>-\frac{1}{2}$, $\frac{1}{2}<\nu<\mu+1$.
\end{proposition}

When specialised to the case of the modified Lommel function of the first kind $\mathbf{L}_\nu(x)$, some of the bounds of Theorems \ref{tiger2} and \ref{tiger1} are sharper than those of \cite{gaunt struve int2} for the integral (\ref{intstruve0}). We record these bounds in the following corollary.

\begin{corollary}\label{struvecor}Let $0<\beta<1$. Then, for $x>0$,
\begin{align}\label{ineqb10ii}\int_0^x\mathrm{e}^{-\beta u}u^\nu \mathbf{L}_\nu(u)\,\mathrm{d}u&<\frac{2\nu+3+\sqrt{2(2\nu+3)}}{(2\nu+1)(1-\beta)}\mathrm{e}^{-\beta x}x^\nu \mathbf{L}_{\nu+1}(x), \quad \nu>-\tfrac{1}{2}, \\
\label{ineqb10iiaa}\int_0^x\mathrm{e}^{-\beta u}u^\nu \mathbf{L}_\nu(u)\,\mathrm{d}u&<\frac{2\nu+3+2\sqrt{2}}{(2\nu+1)(1-\beta)}\mathrm{e}^{-\beta x}x^\nu \mathbf{L}_{\nu+1}(x), \quad \nu>-\tfrac{1}{2}, \\
\int_0^x\mathrm{e}^{-\beta u}u^\nu \mathbf{L}_{\nu}(u)\,\mathrm{d}u&>\frac{1}{1-\beta}\bigg\{\bigg(1-\frac{2\nu(2\nu+1+2\sqrt{2})}{(2\nu-1)(1-\beta)}\frac{1}{x}\bigg)\mathrm{e}^{-\beta x}x^\nu \mathbf{L}_{\nu}(x)\nonumber\\
\label{ineqb23ii}&\quad -\frac{\gamma(2\nu+1,\beta x)}{\sqrt{\pi}2^\nu\beta^{2\nu+1}\Gamma(\nu+\frac{3}{2})}\bigg\}, \quad \nu>\tfrac{1}{2}.
\end{align}
\end{corollary}





\section{Discussion and comparison of the bounds}\label{secnew}

In this section, we discuss the performance of the bounds of Theorems \ref{tiger2} and \ref{tiger1}, and compare our bounds with those given by \cite{gaunt lommel1} for the integral (\ref{intstruve}). We also comment on the improvements the bounds of Corollary \ref{struvecor} made on some of the bounds of \cite{gaunt struve int2} for the integral (\ref{intstruve0}). In this section $0<\beta<1$. 

Subject to the conditions on $\mu$ and $\nu$ given in Theorem \ref{tiger2}, it can be checked that $A_{\mu,\nu}<\frac{1}{2}(3+\sqrt{17})$ (in which case inequality (\ref{ineqb1ss}) is sharper than (\ref{ineqb10})) if $-\frac{3}{2}<\mu\leq\frac{\sqrt{17}}{2}-1$, $-\frac{1}{2}<\nu<\mu+1$, or if $\frac{\sqrt{17}}{2}-1<\mu<\frac{\sqrt{17}}{2}$, $\frac{3+\sqrt{17}}{2}u-\frac{\sqrt{17}+13}{4}<\nu<\mu+1$. In particular, in the important special case $\mu=\nu$ (corresponding to the case of the modified Struve function $\mathbf{L}_\nu(x)$), we have that $A_{\nu,\nu}=\sqrt{2(2\nu+3)}<2\sqrt{2}<\frac{1}{2}(3+\sqrt{17})$ if $-\frac{1}{2}<\nu<\frac{1}{2}$ (recall that $v>-\frac{1}{2}$ is an assumption of Theorem \ref{tiger2}). This fact is helpful in interpreting Corollary \ref{struvecor}. 

Inequality (\ref{ineqb10}) takes a simpler form than (\ref{ineqb1ss}), and so, whilst for some values of $\mu$ and $\nu$ inequality (\ref{ineqb1ss}) is sharper, we consider (\ref{ineqb10}) to be a more preferable general bound than (\ref{ineqb1ss}). Indeed, in deriving inequality (\ref{ineqb12}) in Theorem \ref{tiger1}, we made use of inequality (\ref{ineqb10}). A similar inequality involving $A_{\mu,\nu}$ can be derived by modifying the proof of inequality (\ref{ineqb12}) by instead using (\ref{ineqb1ss}), although for sake of brevity we omit the details.

Inequality (\ref{ineqb4}) is sharper than the only other lower bound for the integral (\ref{intstruve}) of \cite{gaunt lommel1}, $\int_0^x \mathrm{e}^{-\beta u}u^{\nu}\tilde{t}_{\mu,\nu}(u)\,\mathrm{d}u>\mathrm{e}^{-\beta x}x^\nu \tilde{t}_{\mu+1,\nu+1}(x)$, $x>0$, $\mu>-\frac{5}{2}$, $-\mu-2<\nu\leq\mu+3$; this lower bound is the first term in the infinite series of the lower bound (\ref{ineqb4}). The other lower bounds of Theorem \ref{tiger1} all perform worse than (\ref{ineqb4}) and the bound of \cite{gaunt lommel1} for `small' $x$, all being negative for sufficiently small $x$. This is readily seen for (\ref{ineqb3}) and (\ref{ineqb12}). Also, a simple asymptotic analysis of the lower bound (\ref{ineqb2}) using (\ref{Itend0}) shows that, for $-\frac{1}{2}<\nu<0$, the limiting form of this bound is $\frac{\nu}{(\mu+\nu+1)(1-\beta)}\frac{x^{\mu+\nu+1}}{2^\mu\Gamma((\mu-\nu+3)/2)\Gamma((\mu+\nu+3)/2)}<0$, as $x\downarrow0$, whilst, for the case $\nu=0$ the bound is again negative for sufficiently small $x$: the limiting form of the bound is $ -\frac{2\beta (x/2)^{\mu+2}}{(1-\beta)(\mu+2)[\Gamma((\mu+3)/2)]^2}<0$, as $x\downarrow0$.  The bounds (\ref{ineqb2})--(\ref{ineqb12}) have an improved performance for `large' $x$, though. In contrast to the bound of \cite{gaunt lommel1}, these bounds are tight in the limit $x\rightarrow\infty$, and enjoy this property without needing an infinite sum involving modified Lommel functions of the first kind as given in (\ref{ineqb4}).

We now make a comparison to the upper bounds of \cite{gaunt lommel1}. Inequality (2.5) of \cite{gaunt lommel1} gives that, for $x>0$, $\mu>-\frac{1}{2}$, $\frac{1}{2}\leq\nu<\mu+1$,
\begin{align}\int_0^x \mathrm{e}^{-\beta u}u^{\nu}\tilde{t}_{\mu,\nu}(u)\,\mathrm{d}u&<\frac{\mathrm{e}^{-\beta x}x^\nu}{(2\nu+1)(1-\beta)}\bigg(2(\nu+1)\tilde{t}_{\mu+1,\nu+1}(x)-\tilde{t}_{\mu+3,\nu+3}(x) \nonumber\\
&\quad-\frac{x^{\mu+2}}{2^{\mu+2}(\mu+\nu+2)\Gamma\big(\frac{\mu-\nu+1}{2}\big)\Gamma\big(\frac{\mu+\nu+5}{2}\big)}\bigg),\nonumber \\
\label{gau1}&<\frac{2(\nu+1)}{(2\nu+1)(1-\beta)}\mathrm{e}^{-\beta x}x^\nu \tilde{t}_{\mu+1,\nu+1}(x).
\end{align}
Under the same conditions on $\mu$ and $\nu$, another upper bound can be obtained by combining inequalities (2.2) and (2.4) of \cite{gaunt lommel1}: for $x>0$,
\begin{equation}\label{gau2}\int_0^x \mathrm{e}^{-\beta u}u^{\nu}\tilde{t}_{\mu,\nu}(u)\,\mathrm{d}u<\frac{1}{1-\beta}\mathrm{e}^{-\beta x}x^\nu \tilde{t}_{\mu,\nu}(x).
\end{equation}
Inequalities (\ref{ineqb1ss}) and (\ref{ineqb10}) of Theorem \ref{tiger2} increase the range of validity of the bounds (\ref{gau1}) and (\ref{gau2}) to $\mu>-\frac{3}{2}$, $-\frac{1}{2}<\nu<\mu+1$ at the expense of larger multiplicative constants. That the range of validity of inequality (\ref{gau2}) has been increased by inequalities (\ref{ineqb1ss}) and (\ref{ineqb10}) at the cost of larger multiplicative constants can be seen from (\ref{Imon}).
In establishing (\ref{gau1}) and (\ref{gau2}), \cite{gaunt lommel1} used the inequality $\tilde{t}_{\mu,\nu}(x)<\tilde{t}_{\mu-1,\nu-1}(x)$, which holds for $x>0$, $\mu>-\frac{1}{2}$, $\frac{1}{2}\leq\nu<\mu+1$ (see (\ref{Imon})). In the larger parameter regime $\mu>-\frac{3}{2}$, $-\frac{1}{2}<\nu<\mu+1$, this inequality is no longer available to us, and we arrived at the inequalities of Theorem \ref{tiger2} through an alternative approach that bypasses the use of the inequality used by \cite{gaunt lommel1} that results in these larger multiplicative constants. 

When specialised to the case $\mu=\nu$ corresponding to the modified Struve function of the first kind $\mathbf{L}_\nu(x)$ our upper bounds (\ref{ineqb10ii}) (for $|\nu|<\frac{1}{2}$, in which case this is the sharper of the bounds) and (\ref{ineqb10iiaa})  of Corollary \ref{struvecor} do, however, have smaller multiplicative constants than the upper bounds (2.4) and (2.5) of \cite{gaunt struve int2} for the integral (\ref{intstruve0}), with these bounds being $\frac{2\nu+29}{(2\nu+1)(1-\beta)}\mathrm{e}^{-\beta x}x^\nu \mathbf{L}_{\nu+1}(x)$ and $\frac{2\nu+15}{(2\nu+1)(1-\beta)}\mathrm{e}^{-\beta x}x^\nu \mathbf{L}_\nu(x)$, for $x>0$, $\nu>-\frac{1}{2}$.
For sake of comparison with our bounds, note that $\mathbf{L}_{\nu+1}(x)<\mathbf{L}_\nu(x)$, for $x>0$, $\nu>-\frac{1}{2}$ (see \cite{bp14}). 
The improved multiplicative constants are a consequence of a more direct approach given in this paper, whereas the bounds of \cite{gaunt struve int2} were obtained through a series of lemmas that resulted in a build up of errors. In a similar manner, our lower bound (\ref{ineqb23ii}) improves on the lower bound (2.6) of \cite{gaunt struve int2}.
We note that combining (\ref{gau1}) with $\mu=\nu$ (and using that $\tilde{t}_{\nu,\nu}(x)=\mathbf{L}_\nu(x)$) and our bound (\ref{ineqb10iiaa}) gives that, for $x>0$,
\begin{equation*}\int_0^x \mathrm{e}^{-\beta u}u^{\nu}\mathbf{L}_\nu(u)\,\mathrm{d}u<\frac{C_{\nu}}{(2\nu+1)(1-\beta)}\mathrm{e}^{-\beta x}x^\nu \mathbf{L}_\nu(x), \quad \nu>-\tfrac{1}{2},
\end{equation*}
where $C_{\nu}=2(\nu+1)$ for $\nu\geq\frac{1}{2}$, and $C_{\nu}=2\nu+3+\sqrt{2(2\nu+3)}$ for $|\nu|<\frac{1}{2}$. Similar inequalities involving the function $\tilde{t}_{\mu,\nu}(x)$ can be obtained by combining the bounds of Theorem \ref{tiger2} and inequality (\ref{gau1}).


Various two-sided inequalities for the integral (\ref{intstruve}) can be obtained by combining the bounds obtained in this paper and those presented in this section.
As an example, combining (\ref{ineqb4}) and (\ref{gau2}) gives that, for $x>0$ $\mu>-\frac{1}{2}$, $\frac{1}{2}\leq\nu<\mu+1$,
\begin{equation}\label{gau3}\sum_{k=0}^\infty \beta^k \tilde{t}_{\mu+k+1,\nu+k+1}(x)<\frac{\mathrm{e}^{\beta x}}{x^\nu}\int_0^x \mathrm{e}^{-\beta u}u^{\nu}\tilde{t}_{\mu,\nu}(u)\,\mathrm{d}u<\frac{\tilde{t}_{\mu,\nu}(x)}{1-\beta}.
\end{equation}
We used \emph{Mathematica} to calculate the relative error in approximating the normalised integral $F_{\mu,\nu,\beta}(x)=\mathrm{e}^{\beta x}x^{-\nu}\int_0^x \mathrm{e}^{-\beta u}u^{\nu}\tilde{t}_{\mu,\nu}(u)\,\mathrm{d}u$ by the upper bound $U_{\mu,\nu,\beta}(x)$ in (\ref{gau3}), and the lower bound truncated at the fifth term,  $L_{\mu,\nu,\beta}(x)=\sum_{k=0}^4 \beta^k \tilde{t}_{\mu+k+1,\nu+k+1}(x)$. We considered three cases of $\mu-\nu=k$, $k=-0.5,2,5$, and also varied $\nu=1,5,10$ and $\beta=0.25,0.5$. The results are given in Tables \ref{table1} and \ref{table2}. (Tables for the case $\mu=\nu$ are given in \cite{gaunt struve int2}.) For fixed $\mu$, $\nu$ and $x$, we see that increasing $\beta$ from 0.25 to 0.5 increases the relative error in approximating $F_{\mu,\nu,\beta}(x)$ by both $L_{\mu,\nu,\beta}(x)$ and $U_{\mu,\nu,\beta}(x)$. We also see that the bounds are most accurate in the case $\mu-\nu=-0.5$. The lower and upper bounds in (\ref{gau3}) are both tight in the limit $x\rightarrow\infty$, and as expected we observe that, for fixed $\mu$, $\nu$ and $\beta$, the relative error in approximating $F_{\nu,\beta}(x)$ by $U_{\mu,\nu,\beta}(x)$ decreases as $x$ increases. However, the truncated sum $L_{\mu,\nu,\beta}(x)$ is not tight as $x\rightarrow\infty$, with the effect being most pronounced for larger $\beta$. Indeed, using (\ref{eqeq1}) and (\ref{Itendinfinity}), we have that, for $\mu>-\frac{1}{2}$, $\frac{1}{2}\leq\nu<\mu+1$, $\lim_{x\rightarrow\infty}\big(1-\frac{L_{\mu,\nu,0.25}(x)}{F_{\mu,\nu,0.25}(x)}\big)=9.766\times 10^{-4}$ and $\lim_{x\rightarrow\infty}\big(1-\frac{L_{\mu,\nu,0.5}(x)}{F_{\mu,\nu,0.5}(x)}\big)=0.03125$.
 We observe that $U_{\mu,\nu,\beta}(x)$ performs poorly for `small' $x$. This because it is of the wrong asymptotic order as $x\downarrow0$; from (\ref{Itend0}) we have that $\frac{U_{\mu,\nu,\beta}(x)}{F_{\mu,\nu,\beta}(x)}\sim\frac{\mu+\nu+2}{(1-\beta)x}$, as $x\downarrow0$.  The lower bound $L_{\nu,\beta}(x)$ is, however, of the correct asymptotic order as $x\downarrow0$, with $\lim_{x\downarrow0}\big(1-\frac{L_{\mu,\nu,\beta}(x)}{F_{\mu,\nu,\beta}(x)}\big)=\frac{1}{\mu+\nu+3}$, and performs better for `small' $x$.
\begin{table}[h]
\centering
\caption{\footnotesize{Relative error in approximating $F_{\mu,\nu,\beta}(x)$ by $L_{\mu,\nu,\beta}(x)$.}}
\label{table1}
{\scriptsize
\begin{tabular}{|c|rrrrrrr|}
\hline
 \backslashbox{$(\mu,\nu,\beta)$}{$x$}      &    0.5 &    5 &    10 &    15 &    25 &    50 & 100   \\
 \hline
$(0.5,1,0.25)$ & 0.2280 & 0.2066 & 0.1419 & 0.1028 & 0.0656 & 0.0346 & 0.0182 \\
$(4.5,5,0.25)$ & 0.0812 & 0.0853 & 0.0778 & 0.0670 & 0.0503 & 0.0302 & 0.0169  \\
$(9.5,10,0.25)$ & 0.0449 & 0.0474 & 0.0471 & 0.0445 & 0.0378 & 0.0257 & 0.0155 \\   
  \hline
$(3,1,0.25)$ & 0.1461 & 0.1591 & 0.1351 & 0.1024 & 0.0656 & 0.0346 & 0.0182  \\
$(7,5,0.25)$ & 0.0676 & 0.0737 & 0.0737 & 0.0664 & 0.0503 & 0.0302 & 0.0169 \\
$(12,10,0.25)$ & 0.0404 & 0.0431 & 0.0447 & 0.0438 & 0.0378 & 0.0257 & 0.0155 \\ 
\hline
$(6,1,0.25)$ & 0.1019 & 0.1151 & 0.1158 & 0.0991 & 0.0656 & 0.0346 & 0.0182   \\
$(10,5,0.25)$ & 0.0562 & 0.0615 & 0.0650  & 0.0633 & 0.0503 & 0.0302 & 0.0169  \\ 
$(15,10,0.25)$ & 0.0360 & 0.0385 & 0.0406 & 0.0414 & 0.0376 & 0.0257 & 0.0155 \\  
  \hline
   \hline
$(0.5,1,0.5)$ & 0.2348 & 0.2723 & 0.2280 & 0.1845 & 0.1341 & 0.0869 & 0.0602 \\
$(4.5,5,0.5)$ & 0.0825 & 0.1000 & 0.1047 & 0.1005 & 0.0881 & 0.0680 & 0.0522  \\
$(9.5,10,0.5)$ & 0.0453 & 0.0524 & 0.0573 & 0.0591 & 0.0580 & 0.0515 & 0.0440 \\   
  \hline
$(3,1,0.5)$ & 0.1497 & 0.2011 & 0.2096 & 0.1821 & 0.1341 & 0.0869 & 0.0602 \\
$(7,5,0.5)$ & 0.0685 & 0.0849 & 0.0976 & 0.0900 & 0.0881 & 0.0680 & 0.0522  \\
$(12,10,0.5)$ & 0.0407 & 0.0473 & 0.0539 & 0.0578 & 0.0579 & 0.0515 & 0.0440 \\ 
\hline
$(6,1,0.5)$ & 0.1038 & 0.1396 & 0.1696 & 0.1708 & 0.1339 & 0.0869 & 0.0602    \\
$(10,5,0.5)$ & 0.0569 & 0.0696 & 0.0836 & 0.0923 & 0.0879 & 0.0680 & 0.0522  \\ 
$(15,10,0.5)$ & 0.0363 & 0.0418 & 0.0483 & 0.0539 & 0.0576 & 0.0515 & 0.0440 \\  
  \hline
\end{tabular}}
\end{table}

\begin{table}[h]
\centering
\caption{\footnotesize{Relative error in approximating $F_{\mu,\nu,\beta}(x)$ by $U_{\mu,\nu,\beta}(x)$.}}
\label{table2}
{\scriptsize
\begin{tabular}{|c|rrrrrrr|}
\hline
 \backslashbox{$(\mu,\nu,\beta)$}{$x$}      &    0.5 &    5 &    10 &    15 &    25 &    50 & 100   \\
 \hline
$(0.5,1,0.25)$ & 8.1497 & 0.2771 & 0.0872 & 0.0520 & 0.0292 & 0.0139 & 0.0068 \\
$(4.5,5,0.25)$ & 29.3965 & 2.1107 & 0.8520 & 0.5130 & 0.2806 & 0.1300 & 0.0625  \\
$(9.5,10,0.25)$ & 56.0364 & 4.6324 & 1.9741 & 1.1881 & 0.6377 & 0.2868 & 0.1351 \\   
  \hline
$(3,1,0.25)$ & 14.7434 & 0.6111 & 0.1129 & 0.0531 & 0.0292 & 0.0139 & 0.0068 \\
$(7,5,0.25)$ & 36.0379 & 2.5889 & 0.9315 & 0.5206 & 0.2807 & 0.1300 & 0.0625 \\
$(12,10,0.25)$ & 62.6900 & 5.1844 & 2.1181 & 1.2153 & 0.6380 & 0.2868 & 0.1351 \\ 
\hline
$(6,1,0.25)$ & 22.7139 & 1.2479 & 0.2457 & 0.0682 & 0.0292 & 0.0139 & 0.0068   \\
$(10,5,0.25)$ & 44.0269 & 3.3075 & 1.1618 & 0.5736 & 0.2811 & 0.1300 & 0.0625  \\ 
$(15,10,0.25)$ & 70.6841 & 5.9372 & 2.4131 & 1.3226 & 0.6415 & 0.2868 & 0.1351 \\  
  \hline
   \hline
$(0.5,1,0.5)$ & 12.3403 & 0.4845 & 0.1485 & 0.0836 & 0.0452 & 0.0212 & 0.0103 \\
$(4.5,5,0.5)$ & 44.1357 & 3.2057 & 1.3045 & 0.7861 & 0.4286 & 0.1972 & 0.0943  \\
$(9.5,10,0.5)$ & 84.0773 & 6.9721 & 2.9816 & 1.7983 & 0.9664 & 0.4339 & 0.2037 \\   
  \hline
$(3,1,0.5)$ & 22.1891 & 0.9907 & 0.2014 & 0.0879 & 0.0452 & 0.0212 & 0.0103 \\
$(7,5,0.5)$ & 54.0910 & 3.9205 & 1.4276 & 0.7992 & 0.4286 & 0.1972 & 0.0943 \\
$(12,10,0.5)$ & 94.0552 & 7.7986 & 3.1985 & 1.8404 & 0.9667 & 0.4339 & 0.2037 \\ 
\hline
$(6,1,0.5)$ & 34.1224 & 1.9315 & 0.4148 & 0.1205 & 0.0455 & 0.0212 & 0.0103   \\
$(10,5,0.5)$ & 66.0687 & 4.9932 & 1.7741 & 0.8834 & 0.4297 & 0.1972 & 0.0943  \\ 
$(15,10,0.5)$ & 106.0445 & 8.9256 & 3.6404 & 2.0030 & 0.9727 & 0.4339 & 0.2037 \\  
  \hline
\end{tabular}}
\end{table}

\section{Proofs}\label{sec3}

\noindent{\emph{Proof of Lemma \ref{lem1}.}} (i) Let $\mu>-\frac{3}{2}$, $-\frac{1}{2}<\nu<\mu+1$ and $0<\beta<1$. Fix $x_*>0$, and suppose $0<x\leq x_*$. Consider the function
\begin{equation*}V_{\mu,\nu,\beta}(x)=\frac{\mathrm{e}^{\beta x}}{x^\nu \tilde{t}_{\mu+1,\nu+1}(x)}\int_0^x\mathrm{e}^{-\beta u} u^\nu \tilde{t}_{\mu,\nu}(u)\,\mathrm{d}u.
\end{equation*}
Note that the conditions on $\mu$ and $\nu$ ensure that the integral exists. We now observe that
\begin{equation*}\frac{\partial V_{\mu,\nu,\beta}(x)}{\partial \beta}=\frac{\mathrm{e}^{\beta x}}{x^\nu \tilde{t}_{\mu+1,\nu+1}(x)}\int_0^x(x-u)\mathrm{e}^{-\beta u} u^\nu \tilde{t}_{\mu,\nu}(u)\,\mathrm{d}u>0,
\end{equation*}
meaning that $V_{\mu,\nu,\beta}(x)$ is an increasing function of $\beta$. Therefore, for $0<\beta<1$,
\begin{align*}V_{\mu,\nu,\beta}(x)&< \frac{\mathrm{e}^{ x}}{x^\nu \tilde{t}_{\mu+1,\nu+1}(x)}\int_0^{x}\mathrm{e}^{- u} u^\nu \tilde{t}_{\mu,\nu}(u)\,\mathrm{d}u<\frac{x}{2\nu+1}\bigg(\frac{\tilde{t}_{\mu,\nu}(x)}{\tilde{t}_{\mu+1,\nu+1}(x)}+1\bigg) \\
&<\frac{x}{2\nu+1}\bigg(\frac{\mu+\nu+3+x}{x}+1\bigg)=\frac{\mu+\nu+3+2x}{2\nu+1}\leq \frac{\mu+\nu+3+2x_*}{2\nu+1},
\end{align*}
where the second inequality is clear from the integral formula (\ref{intfor}) (which we can apply because $\nu>-\frac{1}{2}$) and we applied inequality (\ref{tilderatio}) to obtain the third inequality (which we can apply because our assumptions on $\mu$ and $\nu$ ensure that $\mu>-2$, $-1\leq\nu<\mu+1$). This completes the proof of inequality (\ref{ineqb1ccc}).

\vspace{1mm}

\noindent{(ii)} Now, fix $x_*>\frac{1}{1-\beta}$. 
We consider the function
\begin{equation*}W_{\mu,\nu,\beta}(x)=M_{\mu,\nu,\beta}(x_*)\mathrm{e}^{-\beta x}x^\nu \tilde{t}_{\mu+1,\nu+1}(x)-\int_0^x \mathrm{e}^{-\beta u}u^\nu \tilde{t}_{\mu,\nu}(u)\,\mathrm{d}u,
\end{equation*}
and prove inequality (\ref{ineqb1}) by showing that $W_{\mu,\nu,\beta}(x)>0$ for all $x\geq x_*$. 

We first prove that $W_{\mu,\nu,\beta}(x_*)>0$, for which it suffices to prove that $V_{\mu,\nu,\beta}(x_*)<M_{\mu,\nu,\beta}(x_*)$. From part (i) and the definition of $M_{\mu,\nu,\beta}(x_*)$ (see (\ref{mng})) we have that
\begin{equation*}V_{\mu,\nu,\beta}(x_*)< \frac{\mu+\nu+3+2x_*}{2\nu+1}\leq M_{\mu,\nu,\beta}(x_*),
\end{equation*}
as required. We now prove that $W_{\mu,\nu,\beta}'(x)>0$ for $x>x_*$.  A calculation using the differentiation formula (\ref{diffone}) followed  by an application of inequality (\ref{Imon}) (we may do so because we assumed $\mu>-\frac{3}{2}$, $-\frac{1}{2}<\nu<\mu+1$) gives that
\begin{align*}W_{\mu,\nu,\beta}'(x)&=M_{\mu,\nu,\beta}(x_*)\frac{\mathrm{d}}{\mathrm{d}x}\big(\mathrm{e}^{-\beta x}x^{-1}\cdot x^{\nu+1} \tilde{t}_{\mu+1,\nu+1}(x)\big)- \mathrm{e}^{-\beta x}x^\nu \tilde{t}_{\mu,\nu}(x)\\
&=M_{\mu,\nu,\beta}(x_*)\mathrm{e}^{-\beta x}x^\nu\big(\tilde{t}_{\mu,\nu}(x)-(x^{-1}+\beta)\tilde{t}_{\mu+1,\nu+1}(x)\big)-\mathrm{e}^{-\beta x}x^\nu \tilde{t}_{\mu,\nu}(x)\\
&>M_{\mu,\nu,\beta}(x_*)\mathrm{e}^{-\beta x}x^\nu\big(1-\beta -x^{-1})\tilde{t}_{\mu,\nu}(x)-\mathrm{e}^{-\beta x}x^\nu \tilde{t}_{\mu,\nu}(x)\\
&\geq\bigg(\frac{1-\beta-x^{-1}}{1-\beta-x_*^{-1}}-1\bigg)\mathrm{e}^{-\beta x}x^\nu \tilde{t}_{\mu,\nu}(x)>0,
\end{align*}
for $x>x_*$. This completes the proof of inequality (\ref{ineqb1}). \hfill $\square$

\vspace{3mm}

\noindent{\emph{Proof of Theorem \ref{tiger2}.}} (i) Suppose $\mu>-\frac{3}{2}$, $-\frac{1}{2}<\nu<\mu+1$ and $0<\beta<1$. Set $x_*=\frac{C}{1-\beta}$ for some $C>1$ that we will choose later. Let $M_{\mu,\nu,\beta}(x_*)$ be defined as in (\ref{mng}). Then, as $\frac{\mu+\nu+3+2x_*}{2\nu+1}\leq M_{\mu,\nu,\beta}(x_*)$, it follows from inequalities (\ref{ineqb1ccc}) and (\ref{ineqb1}) of Lemma \ref{lem1} that, for all $x>0$,
\begin{align}V_{\mu,\nu,\beta}(x)&=\frac{\mathrm{e}^{\beta x}}{x^\nu \tilde{t}_{\mu+1,\nu+1}(x)}\int_0^x\mathrm{e}^{-\beta u} u^\nu \tilde{t}_{\mu,\nu}(u)\,\mathrm{d}u\nonumber\\
&<M_{\mu,\nu,\beta}(x_*)\nonumber\\
&=\max\bigg\{\frac{1}{2\nu+1}\bigg(\mu+\nu+3+\frac{2C}{1-\beta}\bigg),\frac{C}{(C-1)(1-\beta)}\bigg\}\nonumber\\
\label{vfbfb}&<\frac{1}{1-\beta}\max\bigg\{\frac{\mu+\nu+3+2C}{2\nu+1},\frac{C}{C-1}\bigg\}\\
\label{dftg}&=:\frac{1}{1-\beta}\max\{T_1(C),T_2(C)\}.
\end{align}
Solving $T_1(C)=T_2(C)$ gives that $C=\frac{1}{2}A_{\mu,\nu}$, and substituting this choice of $C$ into (\ref{dftg}) gives us inequality (\ref{ineqb1ss}).

\vspace{1mm}

\noindent{(ii)}  We consider the cases $-\frac{1}{2}<\nu<\frac{1}{2}$ and $\nu\geq\frac{1}{2}$ separately. Firstly, suppose that $-\frac{1}{2}<\nu<\frac{1}{2}$. Then, from (\ref{vfbfb}) we get that, for $x>0$,
\begin{align}V_{\mu,\nu,\beta}(x)&<\frac{1}{(2\nu+1)(1-\beta)}\max\bigg\{\mu+\nu+3+2C,\frac{2C}{C-1}\bigg\}\nonumber\\
\label{vrfvse}&<\frac{1}{(2\nu+1)(1-\beta)}\max\bigg\{\mu+\nu+3+2C,\mu+\nu+2+\frac{2C}{C-1}\bigg\},
\end{align}
where we used that $-\frac{1}{2}<\nu<\frac{1}{2}$ in obtaining the first inequality and that $\mu+\nu>-2$ in obtaining the second inequality. We now choose $C>1$ so that $3+2C=2+\frac{2C}{C-1}$, from which we get $C=\frac{1}{4}(3+\sqrt{17})$
. Plugging this choice of $C$ into inequality (\ref{vrfvse}) shows that inequality (\ref{ineqb10}) holds for $-\frac{1}{2}<\nu<\frac{1}{2}$. We deduce that inequality (\ref{ineqb10}) holds for the case $\nu\geq\frac{1}{2}$ from inequality (\ref{gau1}), because our assumption that $\nu<\mu+1$ guarantees that $2\nu+2<\mu+\nu+\frac{1}{2}(9+\sqrt{17})$. This completes the proof of inequality (\ref{ineqb10}).
\hfill $\square$



\vspace{3mm}

\noindent{\emph{Proof of Theorem \ref{tiger1}.}} (i) Let $x>0$ and suppose that $\mu>-1$, $-\mu-1<\nu\leq0$.  By integration by parts and an application of the differentiation formula (\ref{diffone}) we get
\begin{align*}\int_0^x\mathrm{e}^{-\beta u}u^\nu \tilde{t}_{\mu+1,\nu+1}(u)\,\mathrm{d}u&=-\frac{1}{\beta}\mathrm{e}^{-\beta x}x^\nu \tilde{t}_{\mu,\nu}(x)+\frac{1}{\beta}\int_0^x \mathrm{e}^{-\beta u}u^\nu \tilde{t}_{\mu-1,\nu-1}(u)\,\mathrm{d}u,
\end{align*}
where we used the limit $\lim_{x\downarrow0}x^\nu\tilde{t}_{\mu,\nu}(x)=0$, for $\mu+\nu>-1$ (see \ref{Itend0})). The integrals exist for $\mu+\nu>-1$ (see (\ref{Itend0})).  Applying the identity (\ref{Iidentity}) and rearranging now yields
\begin{align}\int_0^x \mathrm{e}^{-\beta u}u^\nu \tilde{t}_{\mu+1,\nu+1}(u)\,\mathrm{d}u+&2\nu\int_0^x\mathrm{e}^{-\beta u}u^{\nu-1}\tilde{t}_{\mu,\nu}(u)\,\mathrm{d}u-\beta\int_0^x\mathrm{e}^{-\beta u}u^\nu \tilde{t}_{\mu,\nu}(u)\,\mathrm{d}u&\nonumber\\
&=\mathrm{e}^{-\beta x}x^\nu \tilde{t}_{\mu,\nu}(x)-\int_0^x\frac{\mathrm{e}^{-\beta u}u^{\mu+\nu}}{2^\mu\Gamma\big(\frac{\mu-\nu+1}{2}\big)\Gamma\big(\frac{\mu+\nu+3}{2}\big)}\,\mathrm{d}u\nonumber\\
\label{55555}&=\mathrm{e}^{-\beta x}x^\nu \tilde{t}_{\mu,\nu}(x)-B_{\mu,\nu,\beta},
\end{align}
where in the final step we evaluated  $\int_0^x\mathrm{e}^{-\beta u}u^{\mu+\nu}\,\mathrm{d}u=\frac{1}{\beta^{\mu+\nu+1}}\gamma(\mu+\nu+1,\beta x)$.
Bounding the first integral in the above display using inequality (\ref{Imon}) and using our assumption that $\nu\leq0$ now yields inequality (\ref{ineqb2}),
as required.

\vspace{1mm}

\noindent{(ii)} Now suppose that $\mu>\frac{1}{2}$, $\frac{3}{2}\leq\nu<\mu+1$. Rearranging (\ref{55555}) gives us
\begin{align}&\int_0^x \mathrm{e}^{-\beta u}u^{\nu}\tilde{t}_{\mu+1,\nu+1}(u)\,\mathrm{d}u-\beta \int_0^x \mathrm{e}^{-\beta u}u^{\nu}\tilde{t}_{\mu,\nu}(u)\,\mathrm{d}u\nonumber\\
\label{1stint}&\quad=\mathrm{e}^{-\beta x}x^\nu \tilde{t}_{\mu,\nu}(x)-2\nu\int_0^x \mathrm{e}^{-\beta u}u^{\nu-1}\tilde{t}_{\mu,\nu}(u)\,\mathrm{d}u-B_{\mu,\nu,\beta}.
\end{align}
We bound the first integral on the left-hand side in (\ref{1stint}) using inequality (\ref{Imon}), and then divide through by $(1-\beta)$ and use inequality (\ref{Imon}) again to get
\begin{align}&\int_0^x \mathrm{e}^{-\beta u}u^{\nu}\tilde{t}_{\mu,\nu}(u)\,\mathrm{d}u\nonumber\\
&\quad>\frac{1}{1-\beta}\bigg\{\mathrm{e}^{-\beta x}x^\nu \tilde{t}_{\mu,\nu}(x)-2\nu\int_0^x \mathrm{e}^{-\beta u}u^{\nu-1}\tilde{t}_{\mu,\nu}(u)\,\mathrm{d}u-B_{\mu,\nu,\beta}\bigg\}\nonumber\\
\label{1stint0}&\quad>\frac{1}{1-\beta}\bigg\{\mathrm{e}^{-\beta x}x^\nu \tilde{t}_{\mu,\nu}(x)-2\nu\int_0^x \mathrm{e}^{-\beta u}u^{\nu-1}\tilde{t}_{\mu-1,\nu-1}(u)\,\mathrm{d}u-B_{\mu,\nu,\beta}\bigg\}.
\end{align}
Finally, we bound the integral $\int_0^x\mathrm{e}^{-\beta u}u^{\nu-1} \tilde{t}_{\mu-1,\nu-1}(u)\,\mathrm{d}u$ using inequality (\ref{gau1}) (we may do so since $\mu>\frac{1}{2}$, $\frac{3}{2}\leq\nu<\mu+1$), which yields inequality (\ref{ineqb3}).

\vspace{1mm}

\noindent{(iii)} We proceed as we did in proving inequality (\ref{ineqb3}), with the only difference being that we use inequality (\ref{ineqb10}) to bound the integral on the right-hand side of (\ref{1stint0}), rather than (\ref{gau1}).  

\vspace{1mm}

\noindent{(iv)} Suppose $\mu>-\frac{5}{2}$, $-\mu-2<\nu\leq\mu+3$, which guarantees that all integrals appearing in this proof of inequality (\ref{ineqb4}) exist and have a positive integrand.  We begin with the same integration by parts as in part (i) of the proof, but with $\mu$ and $\nu$ replaced by $\mu+1$ and $\nu+1$, respectively:
\begin{align}\label{reart}\int_0^x \mathrm{e}^{-\beta u}u^{\nu+1}\tilde{t}_{\mu+1,\nu+1}(u)\,\mathrm{d}u&=-\frac{1}{\beta}\mathrm{e}^{-\beta x}x^{\nu+1}\tilde{t}_{\mu+1,\nu+1}(x)\nonumber\\
&\quad+\frac{1}{\beta}\int_0^x\mathrm{e}^{-\beta u}u^{\nu+1}\tilde{t}_{\mu,\nu}(u)\,\mathrm{d}u.
\end{align}
We now apply the simple inequality $\int_0^x\mathrm{e}^{-\beta u}u^{\nu+1}\tilde{t}_{\mu,\nu}(u)\,\mathrm{d}u<x\int_0^x\mathrm{e}^{-\beta u}u^{\nu}\tilde{t}_{\mu,\nu}(u)\,\mathrm{d}u$, $x>0$, to (\ref{reart}) and rearrange to get
\begin{equation}\label{jj27}\int_0^x\mathrm{e}^{-\beta u}u^{\nu}\tilde{t}_{\mu,\nu}(u)\,\mathrm{d}u>\mathrm{e}^{-\beta x}x^{\nu}\tilde{t}_{\mu+1,\nu+1}(x)+\frac{\beta}{x}\int_0^x \mathrm{e}^{-\beta u}u^{\nu+1}\tilde{t}_{\mu+1,\nu+1}(u)\,\mathrm{d}u.
\end{equation}
From (\ref{jj27}) we can deduce another inequality
\begin{align*}\int_0^x\mathrm{e}^{-\beta u}u^{\nu}\tilde{t}_{\mu,\nu}(u)\,\mathrm{d}u&>\mathrm{e}^{-\beta x}x^{\nu}\tilde{t}_{\mu+1,\nu+1}(x)+\frac{\beta}{x}\bigg(\mathrm{e}^{-\beta x}x^{\nu+1}\tilde{t}_{\mu+2,\nu+2}(x)\\
&\quad+\frac{\beta}{x}\int_0^x \mathrm{e}^{-\beta u}u^{\nu+2}\tilde{t}_{\mu+2,\nu+2}(u)\,\mathrm{d}u\bigg)\\
&=\mathrm{e}^{-\beta x}x^{\nu}\tilde{t}_{\mu+1,\nu+1}(x)+\beta \mathrm{e}^{-\beta x}x^{\nu}\tilde{t}_{\mu+2,\nu+2}(x)\\
&\quad+\frac{\beta^2}{x^2}\int_0^x \mathrm{e}^{-\beta u}u^{\nu+2}\tilde{t}_{\mu+2,\nu+2}(u)\,\mathrm{d}u,
\end{align*}
and iterating yields inequality (\ref{ineqb4}).  In carrying out this iteration, it should be noted that $\sum_{k=0}^\infty \beta^k \tilde{t}_{\mu+k+1,\nu+k+1}(x)$ is a convergent series. To see this, we apply inequality (\ref{Imon}) to obtain that, for all $x>0$, $\sum_{k=0}^\infty \beta^k \tilde{t}_{\mu+k+1,\nu+k+1}(x)<\tilde{t}_{\mu+1,\nu+1}(x)\sum_{k=0}^\infty \beta^k=\frac{\tilde{t}_{\mu+1,\nu+1}(x)}{1-\beta}$, since $0<\beta<1$.



\vspace{1mm}

\noindent{(v)} The tightness of inequalities (\ref{ineqb2})--(\ref{ineqb12}) in the limit $x\rightarrow\infty$ is immediate from the following limiting forms, which hold provided $\mu+\nu>-2$ and $0<\beta<1$:
\begin{align}\label{eqeq1} \int_0^x \mathrm{e}^{-\beta u}u^\nu  \tilde{t}_{\mu,\nu}(u)\,\mathrm{d}u&\sim \frac{1}{\sqrt{2\pi}(1-\beta)}x^{\nu-1/2}\mathrm{e}^{(1-\beta)x}, \quad x\rightarrow\infty,\\
\label{eqeq2}\mathrm{e}^{-\beta x}x^\nu \tilde{t}_{\mu+n,\nu+n}(x)&\sim  \frac{1}{\sqrt{2\pi}}x^{\nu-1/2}\mathrm{e}^{(1-\beta)x}, \quad x\rightarrow\infty,\:n\in\mathbb{R}.
\end{align}
Here, (\ref{eqeq1}) follows from using (\ref{Itendinfinity}) and a standard asymptotic analysis, whilst (\ref{eqeq2}) is immediate from (\ref{Itendinfinity}).  To verify that (\ref{ineqb4}) is tight in the limit $x\rightarrow\infty$ we additionally use that $\sum_{k=0}^\infty\beta^k=\frac{1}{1-\beta}$, since $0<\beta<1$. \hfill $\square$

\vspace{3mm}

\noindent{\emph{Proof of Proposition \ref{propone}.}} (i) 
In part (i) of the proof of Theorem \ref{tiger1}, use inequality (\ref{Imon}) to bound the third integral in (\ref{55555}), rather than the first integral.

\vspace{1mm}

\noindent{(ii)} 
In part (ii) of the proof of Theorem \ref{tiger1}, use (\ref{Imon}) to bound the second integral in (\ref{1stint}), rather than the first integral.

\vspace{1mm}

\noindent{(iii)} On examining the proof of inequality (\ref{ineqb12}), we see that the modification given in part (ii) that gave us inequality (\ref{ineqb22}) rather than (\ref{ineqb3}) can also be used to yield inequality (\ref{ineqb23}). \hfill $\square$

\vspace{3mm}

\noindent{\emph{Proof of Corollary \ref{struvecor}.}} (i) Inequality (\ref{ineqb10ii}) follows from setting $\mu=\nu$ in (\ref{ineqb1ss}) and using that $\tilde{t}_{\nu,\nu}(x)=\mathbf{L}_\nu(x)$. 

\vspace{1mm}

\noindent{(ii)} We prove inequality (\ref{ineqb10iiaa}) by considering the cases $<-\frac{1}{2}<\nu<\frac{1}{2}$ and $\nu\geq\frac{1}{2}$ separately. That inequality (\ref{ineqb10iiaa}) holds for $<-\frac{1}{2}<\nu<\frac{1}{2}$ follows from inequality (\ref{ineqb10ii}) because $\sqrt{2(2\nu+3)}<2\sqrt{2}$ for $<-\frac{1}{2}<\nu<\frac{1}{2}$. We deduce that (\ref{ineqb10iiaa}) holds for $\nu\geq\frac{1}{2}$ from inequality (\ref{gau1}), since $2\nu+2<2\nu+3+2\sqrt{2}$. 

\vspace{1mm}

\noindent{(iii)} Follow the proof of inequality (\ref{ineqb12}) with $\mu=\nu$ (so that $\tilde{t}_{\nu,\nu}(x)=\mathbf{L}_\nu(x)$), but bound the integral on the right-hand side of (\ref{1stint0}) using inequality (\ref{ineqb10iiaa}) instead of inequality (\ref{ineqb10}) (in which we would have set $\mu=\nu$).
\hfill $\square$

\appendix

\section{Basic properties of the modified Lommel function of the first kind}\label{appa}

This appendix contains a list of some basic properties of the modified Lommel function of the first kind $\tilde{t}_{\mu,\nu}(x)$ that we use in this paper. We also define the Kamp\'e de F\'eriet function. The modified Lommel function $\tilde{t}_{\mu,\nu}(x)$ is a regular function of $x\in\mathbb{R}$. For $x>0$, the function $\tilde{t}_{\mu,\nu}(x)$ is positive if $\mu-\nu\geq-3$ and $\mu+\nu\geq-3$. 
The following recurrence relation and differentiation formula for  can be found in \cite{gaunt lommel} and \cite{r64}, respectively:
\begin{align}\label{Iidentity}\tilde{t}_{\mu-1,\nu-1}(x)-\tilde{t}_{\mu+1,\nu+1}(x)&=\frac{2\nu}{x}\tilde{t}_{\mu,\nu}(x)+a_{\mu,\nu}(x), \\
\label{diffone}\frac{\mathrm{d}}{\mathrm{d}x} \big(x^{\nu}  \tilde{t}_{\mu,\nu} (x) \big) &= x^{\nu} \tilde{t}_{\mu-1,\nu -1} (x),
\end{align}
where 
\[a_{\mu,\nu}(x)=\frac{(\frac{1}{2}x)^\mu}{\Gamma\big(\frac{\mu-\nu+1}{2}\big)\Gamma\big(\frac{\mu+\nu+3}{2}\big)}.\] 
The following asymptotic properties are also given in \cite{gaunt lommel}:
\begin{align}\label{Itend0} \tilde{t}_{\mu,\nu}(x)&\sim\frac{(\frac{1}{2}x)^{\mu+1}}{\Gamma\big(\frac{\mu-\nu+3}{2}\big)\Gamma\big(\frac{\mu+\nu+3}{2}\big)}\bigg(1+\frac{x^2}{(\mu+3)^2-\nu^2}\bigg), \quad x\downarrow0, \\
\label{Itendinfinity} \tilde{t}_{\mu,\nu}(x)&\sim\frac{\mathrm{e}^x}{\sqrt{2\pi x}}, \quad x\rightarrow\infty, 
\end{align}
where (\ref{Itend0}) is valid for $\mu>-3$, $|\nu|<\mu+3$, whilst (\ref{Itendinfinity}) holds for all $\mu,\nu\in\mathbb{R}$.

Let $x>0$, $\mu>-\frac{1}{2}$ and $\frac{1}{2}\leq\nu<\mu+1$. Then it was shown by \cite{gaunt lommel} that
\begin{equation}\label{Imon}\tilde{t}_{\mu,\nu}(x)<\tilde{t}_{\mu-1,\nu-1}(x).  
\end{equation} 
Inequality (\ref{Imon}) generalises an inequality of \cite{bp14} for the modified Struve function of the first kind $\mathbf{L}_\nu(x)$. For further functional inequalities involving the modified Lommel function of the first kind, some of which improve on (\ref{Imon}), see \cite{gaunt lommel,gaunt lommel func ineq,m17}. Suppose now that $\mu>-1$, $0\leq\nu<\mu+1$, $x>0$. Then, we also have the following lower bound of \cite[inequality (3.42)]{gaunt lommel},
\begin{equation*}\frac{\tilde{t}_{\mu,\nu}(x)}{\tilde{t}_{\mu-1,\nu-1}(x)}>\frac{x}{\mu+\frac{1}{2}+\sqrt{(\nu+\frac{1}{2})^2+x^2}},
\end{equation*}
from which the following simpler inequality follows by the triangle inequality
\begin{equation}\label{tilderatio}\frac{\tilde{t}_{\mu,\nu}(x)}{\tilde{t}_{\mu-1,\nu-1}(x)}>\frac{x}{\mu+\nu+1+x}.
\end{equation}

The Kamp\'e de F\'eriet function \cite{k37}, written according to the notation of \cite[p.\ 423, equation (26)]{lond}, is defined by
\begin{align}\label{kamp} \mathrm{F}^{H:A;B}_{G:C;D} \Bigg[ \begin{array}{c} (h_H):(a_A);\,(b_B)\\ (g_G):(c_C);\, (d_D)
   \end{array} \Bigg| x, y \Bigg]=\sum_{m=0}^\infty\sum_{n=0}^\infty\frac{((h_H))_{m+n}((a_A))_m((b_B))_n}{((g_G))_{m+n}((c_C))_m((d_D))_n}\frac{x^m}{m!}\frac{y^n}{n!}\,, 
   \end{align}
where $(h_H)$ denotes the sequence $(h_1, h_2, \ldots , h_H )$ and $((h_H ))_n$ is defined by
the product of Pochhammer symbols $((h_H))_n:=(h_1)_n(h_2)_n\cdots(h_H)_n$, $n\geq0$. We use the convention that $((h_H))_0=1$. For further details of the function (\ref{kamp}), see \cite[pp.\ 26--33]{sii}.

\section{Proof of the integral formula (\ref{intfor})}\label{appb}

\noindent{We} wish to prove that $I_{\mu,\nu}(x):=\int_0^x \mathrm{e}^{-u}u^\nu \tilde{t}_{\mu,\nu}(u)\,\mathrm{d}u=f_{\mu,\nu}(x)$,
where
\begin{equation*}f_{\mu,\nu}(x)=\frac{\mathrm{e}^{-x}x^{\nu+1}}{2\nu+1}\big(\tilde{t}_{\mu,\nu}(x)+\tilde{t}_{\mu+1,\nu+1}(x)\big)-\frac{\gamma(\mu+\nu+2,x)}{2^\mu(2\nu+1)\Gamma\big(\frac{\mu-\nu+1}{2}\big)\Gamma\big(\frac{\mu+\nu+3}{2}\big)}.
\end{equation*}
We shall do so by showing that $f_{\mu,\nu}'(x)=\mathrm{e}^{-x}x^\nu \tilde{t}_{\mu,\nu}(x)$, which suffices because $\lim_{x\downarrow0}I_{\mu,\nu}(x)=\lim_{x\downarrow0}f_{\mu,\nu}(x)=0$. We begin by noting that an application of the differentiation formula (\ref{diffone}) gives
\begin{equation*}\frac{\mathrm{d}}{\mathrm{d}x}\big(x^{\nu+1}\tilde{t}_{\mu,\nu}(x)\big)=\frac{\mathrm{d}}{\mathrm{d}x}\big(x\cdot x^{\nu}\tilde{t}_{\mu,\nu}(x)\big)=x^\nu\tilde{t}_{\mu,\nu}(x)+x^{\nu+1}\tilde{t}_{\mu-1,\nu-1}(x).
\end{equation*}
Using this differentiation formula as well as (\ref{Iidentity}), followed by an application of the identity (\ref{diffone}) gives
\begin{align*}(2\nu+1)f_{\mu,\nu}'(x)&=\mathrm{e}^{-x}\Big[-x^{\nu+1}\big(\tilde{t}_{\mu,\nu}(x)+\tilde{t}_{\mu+1,\nu+1}(x)\big)+x^\nu\tilde{t}_{\mu,\nu}(x)\\
&\quad+x^{\nu+1}\tilde{t}_{\mu-1,\nu-1}(x)+x^{\nu+1}\tilde{t}_{\mu,\nu}(x)\Big]-\frac{x^{\mu+\nu+1}\mathrm{e}^{-x}}{2^\mu\Gamma\big(\frac{\mu-\nu+1}{2}\big)\Gamma\big(\frac{\mu+\nu+3}{2}\big)}\\
&=\mathrm{e}^{-x}x^\nu\big[\tilde{t}_{\mu,\nu}(x)-x\tilde{t}_{\mu+1,\nu+1}(x)+x\tilde{t}_{\mu-1,\nu-1}(x)-xa_{\mu,\nu}(x)\big]\\
&=\mathrm{e}^{-x}x^\nu[\tilde{t}_{\mu,\nu}(x)+2\nu\tilde{t}_{\mu,\nu}(x)]=(2\nu+1)\mathrm{e}^{-x}x^\nu\tilde{t}_{\mu,\nu}(x),
\end{align*}
as required. \hfill $\Box$

\section*{Acknowledgements}
The author is supported by a Dame Kathleen Ollerenshaw Research
Fellowship. I would like to thank the reviewer for their helpful comments and suggestions.

\end{document}